\documentclass[12pt,twoside]{article}
\usepackage[latin1]{inputenc}
\usepackage{times}
\usepackage{epsfig}
\usepackage{amsfonts}
\usepackage{float}
\usepackage{amsmath}
\usepackage{latexsym}
\usepackage{graphicx}
\usepackage{nicefrac} 
\usepackage{type1cm}         
\usepackage{txfonts}
\usepackage{xcolor}

\relax

\newcommand{\liff}{\Longleftrightarrow}
\newcommand{\dl}{\displaystyle}

\newcommand{\fts}{\footnotesize}

\newcommand{\RA}{\Rightarrow}

\newcommand{\nnr}{\nonumber}

\newcommand{\sep}{\quad}

\newcommand{\mc}{\multicolumn}

\newcommand{\bay}{\begin{array}}
\newcommand{\eay}{\end{array}}

\newcommand{\bsl}{\begin{slide}}
\newcommand{\esl}{\end{slide}}

\newcommand{\bol}{\begin{overlay}}
\newcommand{\eol}{\end{overlay}}

\newcommand{\bqa}{\begin{eqnarray*}}
\newcommand{\eqa}{\end{eqnarray*}}

\newcommand{\bqan}{\begin{eqnarray}}
\newcommand{\eqan}{\end{eqnarray}}

\newcommand{\bqt}{\begin{quotation}}
\newcommand{\eqt}{\end{quotation}}

\newcommand{\bt}{\begin{tabbing}}
\newcommand{\et}{\end{tabbing}}

\newcommand{\bit}{\begin{itemize}}
\newcommand{\eit}{\end{itemize}}

\newcommand{\bist}{\begin{itemstep}}
\newcommand{\eist}{\end{itemstep}}

\newcommand{\ben}{\begin{enumerate}}
\newcommand{\een}{\end{enumerate}}

\newcommand{\beq}{\begin{equation}}
\newcommand{\eeq}{\end{equation}}

\newcommand{\bdes}{\begin{description}}
\newcommand{\edes}{\end{description}}

\newcommand{\btb}{\begin{tabular}}
\newcommand{\etb}{\end{tabular}}

\newcommand{\bpic}{\begin{picture}}
\newcommand{\epic}{\end{picture}}

\newcommand{\bcen}{\begin{center}}
\newcommand{\ecen}{\end{center}}

\newcommand{\bfg}{\begin{figure}}
\newcommand{\efg}{\end{figure}}

\newcommand{\bmp}{\begin{minipage}}
\newcommand{\emp}{\end{minipage}}

\newcommand{\bgan}{\begin{gather}}
\newcommand{\egan}{\end{gather}}

\newcommand{\bal}{\begin{align*}}
\newcommand{\eal}{\end{align*}}

\newcommand{\baln}{\begin{align}}
\newcommand{\ealn}{\end{align}}

\newcommand{\bala}{\begin{alignat*}}
\newcommand{\eala}{\end{alignat*}}

\newcommand{\balan}{\begin{alignat}}
\newcommand{\ealan}{\end{alignat}}

\newcommand{\bspt}{\begin{split}}
\newcommand{\espt}{\end{split}}

\newcommand{\bvt}{\begin{verbatim}}
\newcommand{\evt}{\end{verbatim}}

\newcommand{\bgb}{\begin{svgraybox}}
\newcommand{\egb}{\end{svgraybox}}

\newcommand{\bdefi}{\begin{definition}[section]}
\newcommand{\edefi}{\end{definition}}

\newcommand{\bco}{\begin{corollary}[definition]}
\newcommand{\eco}{\end{corollary}}

\newcommand{\bpro}{\begin{proposition}[definition]}
\newcommand{\epro}{\end{proposition}}

\newcommand{\blem}{\begin{lemma}[definition]}
\newcommand{\elem}{\end{lemma}}

\newcommand{\bas}{\begin{assumption}[definition]}
\newcommand{\eas}{\end{assumption}}

\newcommand{\babb}{\begin{figure}[section]}
\newcommand{\eabb}{\end{figure}}

\newcommand{\btab}{\begin{table}}
\newcommand{\etab}{\end{table}}

\newcommand{\bexe}{\begin{exercise}}
\newcommand{\eexe}{\end{exercise}}

\newcommand{\bprf}{\begin{proof}}
\newcommand{\eprf}{\qed \end{proof}}

\newcommand{\bprob}{\begin{problem}}
\newcommand{\eprob}{\end{problem}}

\newtheorem{ass}{{\bf Assumptions}}
\newcommand{\bass}{ \begin{ass} }
\newcommand{\eass}{\end{ass}}

\newtheorem{nota}{{\bf Notation}}
\newcommand{\bnot}{ \begin{nota} }
\newcommand{\enot}{\end{nota}}

\newtheorem{bedi}{{\bf Condition}}
\newcommand{\bbd}{ \begin{bedi}}
\newcommand{\ebd}{\end{bedi}}

\newtheorem{remark}{{\sc Remark}}
\newcommand{\brem}{\begin{remark}}
\newcommand{\erem}{\end{remark}}

\newtheorem{notas}{{\bf Notations}}
\newcommand{\bnots}{ \begin{notas}}
\newcommand{\enots}{\end{notas}}

\newtheorem{appropr}{{\bf Approximation Procedure}}
\newcommand{\bappr}{\begin{appropr}}
\newcommand{\eappr}{\end{appropr}}

\newtheorem{result}{{\bf Result}}
\newcommand{\bres}{ \begin{result} }
\newcommand{\eres}{\end{result}}

\newtheorem{lemma}{{\bf Lemma}}

\newtheorem{schema}{{\bf Schematic}}{\rm }
\newcommand{\bscm}{ \begin{schema} }
\newcommand{\escm}{\end{schema}}

\newtheorem{model}{{\bf Model}}{\rm }
\newcommand{\bmdl}{ \begin{model} }
\newcommand{\emdl}{\end{model}}

{\rm }
\newcommand{\bexa}{ \begin{example} }
\newcommand{\eexa}{\end{example}}

\newcommand{\df}{distribution function}
\newcommand{\db}{distribution}
\newcommand{\dfs}{distribution functions}
\newcommand{\dbs}{dis\-tri\-bu\-tions}
\newcommand{\ind}{independent}

\newcommand{\yp}{hypothesis}
\newcommand{\yps}{hypotheses}
\newcommand{\np}{non\-pa\-ra\-me\-tric}
\newcommand{\stat}{statistic}
\newcommand{\stats}{statistics}
\newcommand{\obs}{observation}
\newcommand{\obss}{observations}

\newcommand{\cis}{confidence intervals}

\newcommand{\asy}{asymptotic}
\renewcommand{\ind}{independent}

\newcommand{\cm}{covariance matrix}

\newcommand{\nfrac}{\nicefrac}

\newcommand{\ig}{i=1, \ldots,}

\newcommand{\kg}{k=1, \ldots,}
\newcommand{\rg}{r=1, \ldots,}

\newcommand{\ellg}{\ell=1, \ldots,}

\newcommand{\appeq}{\mbox{ \d{$\stackrel{\textstyle .}{=}$} }}

\newcommand{\sumi}{\sum_{i=1}^}
\newcommand{\sumj}{\sum_{j=1}^}
\newcommand{\sumk}{\sum_{k=1}^}
\newcommand{\sumell}{\sum_{\ell=1}^}

\newcommand{\sumr}{\sum_{r=1}^}

\newcommand{\Cov}{\operatorname{{\it Cov}}}

\newcommand{\hf}{{\widehat f}}
\newcommand{\hF}{{\widehat F}}

\newcommand{\hG}{{\widehat G}}

\newcommand{\hH}{{\widehat H}}

\newcommand{\hp}{{\widehat p}}

\newcommand{\hpsi}{{\widehat \psi}}

\newcommand{\hsigma}{{\widehat \sigma}}

\newcommand{\olp}{\overline{p}}

\newcommand{\olR}{\overline{R}}

\newcommand{\vc}{\boldsymbol{c}}

\newcommand{\vn}{\boldsymbol{n}}

\newcommand{\vp}{\boldsymbol{p}}

\newcommand{\vC}{\boldsymbol{C}}

\newcommand{\vF}{\boldsymbol{F}}

\newcommand{\vI}{\boldsymbol{I}}
\newcommand{\vJ}{\boldsymbol{J}}

\newcommand{\vT}{\boldsymbol{T}}

\newcommand{\vW}{\boldsymbol{W}}

\newcommand{\vmu}{\boldsymbol{\mu}}

\newcommand{\vSigma}{\boldsymbol{\Sigma}}

\newcommand{\vpsi}{\boldsymbol{\psi}}

\newcommand{\veins}{{\bf 1}}
\newcommand{\vnull}{{\bf 0}}

\newcommand{\volR}{\boldsymbol{\overline{R}}}

\newcommand{\volY}{\boldsymbol{\overline{Y}}}
\newcommand{\volZ}{\boldsymbol{\overline{Z}}}

\newcommand{\vwhp}{\boldsymbol{\widehat{p}}}

\newcommand{\vwhF}{\boldsymbol{\widehat{F}}}

\newcommand{\vwhpsi}{\boldsymbol{\widehat{\psi}}}

\begin{document}

\thispagestyle{empty}
\begin{center}
{{\Large \bf Ranks and Pseudo-Ranks} \\[1ex] - {\large \bf Paradoxical Results 
  of Rank Tests -}}
\end{center}
\bcen
{\large Edgar Brunner} \\[1ex]
{\it Department of Medical Statistics, University of G{\"o}ttingen, Germany} 
\\[2ex]
{\large Frank Konietschke}\\[1ex]
{\it Department of Mathematical Sciences, University of Texas at Dallas, 
U.S.A.} \\[2ex]
{\large Arne C. Bathke} \\[1ex]
{\it Department of Mathematics, University of Salzburg, Austria} 
\\[2ex]
and \\[2ex]
{\large Markus Pauly} \\[1ex]
{\it Institute of Statistics, Ulm University, Germany}
\ecen
\mbox{ }\\
\bcen
{\sc Abstract} \\[3mm]
\ecen
Rank-based inference methods are applied in various disciplines, typically when 
procedures relying on standard normal theory are not justifiable, for example 
when data are not symmetrically distributed, contain outliers, or responses are 
even measured on ordinal scales. Various specific rank-based methods have been 
developed for two and more samples, and also for general factorial designs 
(e.g., Kruskal-Wallis test, Jonckheere-Terpstra test). It is the aim of the 
present paper (1) to demonstrate that traditional rank-procedures for several 
samples or general factorial designs may lead to paradoxical results in case of 
unbalanced samples, (2) to explain why this is the case, and (3) to provide a 
way to overcome these disadvantages of traditional rank-based inference. 
Theoretical investigations show that the paradoxical results can be explained 
by carefully considering the non-centralities of the test statistics which may 
be non-zero for the traditional tests in unbalanced designs. These 
non-centralities may even become arbitrarily large for increasing sample sizes 
in the unbalanced case. A simple solution is the use of so-called pseudo-ranks 
instead of ranks. As a special case, we illustrate the effects in sub-group 
analyses which are often used when dealing with rare diseases. \\[4ex]
\newpage
\mbox{ }

\renewcommand{\baselinestretch}{1}




\section{Introduction} \label{int}

If the assumptions of classical parametric inference methods are not met, the 
usual recommendation is to apply nonparametric rank-based tests. 
Here, the Wilcoxon-Mann-Whitney and Kruskal-Wallis (1952) tests are among 
the most commonly applied rank procedures,
often utilized as replacements for the unpaired two-sample $t$-test and 
the one-way ANOVA, respectively. 
Other popular rank methods include the 
Hettmansperger-Norton (1987) and Jonckheere-Terpstra (1952, 1954) tests for 
ordered alternatives, and the procedures by Akritas et al. (1997) for two- or 
higher-way designs. In statistical practice, these procedures are usually 
appreciated as robust and powerful inference tools when standard assumptions 
are not fulfilled. For example, Whitley and Ball (2002) conclude that 
{\it ``Nonparametric methods require no or very limited assumptions to be made 
about the format of the data, and they may, therefore, be preferable when the 
assumptions required for parametric methods are not valid.''} 
In line with this statement, Bewick et al. (2004) also state that 
{\it ``the Kruskal-Wallis, Jonckheere-Terpstra (..) tests can be used to test 
for differences between more than two groups or treatments when the 
assumptions for analysis of variance are not held.''} 

These descriptions are slightly over-optimistic since nonparametric methods also 
rely on certain assumptions. In particular, the Wilcoxon-Mann-Whitney and 
Kruskal-Wallis tests postulate homoscedasticity across groups under the null
hypothesis, and they have originally only been developed for continuous 
outcomes. In case of doubt, it is nevertheless expected that rank procedures 
are more robust and lead to more reliable results than their parametric 
counterparts. While this is true for deviations from normality, and while by 
now it is widely accepted that ordinal data should rather be analyzed using 
adequate rank-based methods than using normal theory procedures, we illustrate 
in various instances that nonparametric rank tests for more than two samples 
possess one noteworthy weakness. Namely, they are generally non-robust against 
changes from balanced to unbalanced designs. In particular, keeping the data 
generating processes fixed, we provide paradigms under which commonly used rank 
tests surprisingly yield {\it completely opposite} test decisions when 
rearranging group sample sizes. These examples are in general not artificially 
generated to obtain paradoxical results, but even include homoscedastic normal 
models. This effect is completely undesirable, leading to the somewhat 
heretical question
\bit
\item Are nonparametric rank procedures useful at all to handle questions for 
more than two groups?
\eit 

In order to comprehensively answer this question, we carefully analyze the 
underlying nonparametric effects of the respective rank procedures. From this, 
we develop detailed guidelines for an adequate application of rank-based 
procedures. Moreover, we even state a simple solution for all these problems: 
Substituting ranks by closely related quantities, the so-called 
{\it pseudo-ranks} that have already been considered by Kulle (1999), Gao and 
Alvo (2005$a,b$), Gao, Alvo, Chen, and Li (2008), and in more detail by 
Thangavelu and Brunner (2007), and by Brunner, Konietschke, Pauly and Puri 
(2017). It should be noted that the motivation in these references was 
different, and that their authors had not been aware of the striking 
paradoxical properties that may arise when using classical rank tests. These 
surprising paradigms only appear in case of unbalanced designs since all rank 
procedures discussed below coincide with their respective pseudo-rank analogs 
in case of equal sample sizes. Pseudo-ranks are easy to compute, share the same 
advantageous properties of ranks and lead to reliable and robust inference 
procedures for a variety of factorial designs. Moreover, we can even obtain 
confidence intervals for (contrasts of) easy to interpret reasonable 
nonparametric effects. Thus, resolving the commonly raised disadvantage that 
{\it ``nonparametric methods are geared toward hypothesis testing rather than 
estimation of effects''} (Whitley and Ball, 2002).

The paper is organized as follows. Notations are introduced in 
Section~\ref{statmod}. Then in Section~\ref{owl} some paradoxical results are 
presented in the one-way layout for the Kruskal-Wallis test and for the 
Hettmansperger-Norton trend test by means of certain tricky (non-transitive) 
dice. In the two-way layout, a paradoxical result for the Akritas-Arnold-Brunner 
test in a simple $2\times 2$-design is presented in Section~\ref{twl} using a 
homoscedastic normal shift model. The theoretical background of the paradoxical 
results is discussed in Section~\ref{expl} and a solution of the problem by 
using pseudo-ranks is investigated in detail. Moreover, the computation of 
confidence intervals is discussed and applied to the data in Section~\ref{twl}. 
Section~\ref{cauno.sg} provides a cautionary note for the problem of 
sub-group analysis where typically unequal sample sizes appear. 

The paper closes with some guidelines for adequate application of rank 
procedures in the discussion and conclusions section. There, it is also briefly 
discussed that the use of pairwise and stratified rankings would make matters 
potentially worse.

\section{Statistical Model and Notations} \label{statmod}

For $d > 2$ samples of $N = \sumi d n_i$ \ind\ \obss\ $X_{ik} \sim F_i = 
\frac12 [F_i^- + F_i^+], \ \ig d, \ \kg n_i$, the \np\ relative effects 
which are underlying the rank tests are commonly defined as  
\bqan
p_i = \int H dF_i, & \text{where} & H = \frac1N \sumr d n_r F_r 
\label{statmod.pi}
\eqan
denotes the weighted mean \db\ of the \dbs\ $F_1, \ldots, F_d$ in the design. 
Here we use the so-called normalized version of the the \db\ $F_i$ to cover the 
cases of continuous, as well as non-continuous \dbs\ in a unified approach. 
Thus, the case of ties does not require a separate consideration.
This idea was first mentioned by Kruskal (1952) and later considered in more 
detail by Ruymgaart (1980). Akritas, Arnold, and Brunner (1997) extended this 
approach to factorial designs, while Akritas and Brunner (1997) and Brunner, 
Munzel, and Puri (1999) applied this technique to repeated measures and 
longitudinal data.

Easily interpreted, $p_i = P(Z < X_{i1}) + \frac12 P(Z = X_{i1})$ is the 
probability that a randomly selected \obs\ $Z$ from the weighted mean \db\ $H$ 
is smaller than a randomly selected \obs\ $X_{i1}$ from the \db\ $F_i$ plus 
$\frac12$ times the probability that both observations are equal. 
Thus, the quantity $p_i$ measures an effect of the \db\ $F_i$ 
with respect to the weighted mean \db\ $H$. 
In the case of two \ind\ random variables $X_1 \sim F_1$ and $X_2 \sim 
F_2$, Birnbaum and Klose (1957) had called the function $L(t) = 
F_2[F_1^{-1}(t)]$ the ``relative \df'' of $X_1$ and $X_2$, 
assuming continuous \dbs. Thus, its expectation
\bqa
 \int_0^1 t d L(t) &=& \int_{- \infty}^\infty F_1(s) d F_2(s) \ = \ 
                       P(X_1 < X_2) 
\eqa
is called a ``relative effect'' with an obvious adaption of the notation. In 
the same way, the quantity $p_i = P(Z < X_{i1}) + \frac12 P(Z = X_{i1})$ is 
called a ``relative effect'' of $X_{i1} \sim F_i$ with respect to 
the weighted mean $Z \sim H$. This effect $p_i$ is a linear combination of the 
pairwise effects $w_{ri} = \int F_r d F_i$. In vector notation, 
Equation~(\ref{statmod.pi}) is written as 
\bqan
 \vp &=& \int H d \vF \ = \ \vW' \vn \ = \ \left( \bay{ccc} w_{11}, & \cdots & 
         w_{d1} \\ \vdots & \ddots & \vdots \\ w_{1d}, & \cdots & w_{dd} \eay
				\right) \cdot \left( \bay{c} n_1/N \\ \vdots \\ n_d/N \eay \right) \ = 
				\ \left( \bay{c} p_1 \\ \vdots \\ p_d \eay \right) \ . \label{vpdef}
\eqan

Here, $\vF = (F_1, \ldots, F_d)'$ denotes the vector of distribution functions, and
\bqan
 \vW &=& \int \vF'\ d \vF \ = \ \left( \bay{ccc} w_{11}, & \cdots & 
         w_{1d} \\ \vdots & \ddots & \vdots \\ w_{d1}, & \cdots & w_{dd} \eay
				\right)
\eqan
is the matrix of pairwise effects $w_{ri}$. Note that $w_{ii} = \frac12$ 
and $w_{ir} = 1 - w_{ri}$ which follows from integration by parts. The relative 
effects $p_i$ are arranged in the vector $\vp = (p_1, \ldots, p_d)'$ and can be 
estimated consistently by the simple plug-in estimator
\bqan 
 \hp_i &=& \int \hH d \hF_i \ = \ \frac1N \left(\olR_{i \cdot} - \tfrac12  
           \right).  \label{piest}
\eqan

Here, $\hF_i$ denotes the (normalized) empirical \db\ of $X_{i1}, \ldots, 
X_{in_i}$, $\ig d$, and $\hH = \frac1N \sumr d n_r \hF_r$ their weighted mean. 
Finally, $\olR_{i \cdot} = \frac1{n_i} \sumk {n_i} R_{ik}$ denotes the mean of 
the ranks 
\bqan
 R_{ik} &=& \frac12 + N \hH(X_{ik}) \ = \ \frac12 + \sumr d \sumell {n_r} 
            c(X_{ik} - X_{r \ell}), \label{rikdef}
\eqan
where the function $c(u) = 0, \nfrac12, 1$ for $u <, =$ or $>0$, respectively,
denotes the count function. The estimators $\hp_1, \ldots, \hp_d$ are arranged 
in the vector 
\bqan
 \vwhp &=& \int \hH d \vwhF \ = \ \frac1N \left( \volR_\cdot - \tfrac12 
           \veins_d \right),                        \label{vechp}
\eqan
where $\vwhF = (\hF_1, \ldots, \hF_d)' $ is the vector of the empirical \dbs, 
$\volR_\cdot = (\olR_{1 \cdot}, \ldots, \olR_{d \cdot})'$ the vector of the 
rank means $\olR_{i \cdot}$, and $\veins_d = (1, \ldots, 1)'_{d\times 1}$ 
denotes the vector of 1s.

In the following sections we demonstrate that for $d \ge 3$ groups, 
rank tests may lead to paradoxical results in case of unequal sample sizes. 
In particular, for factorial designs involving two or more factors, 
the \np\ main effects and interactions (defined by the weighted relative effects 
$p_{ij} = \int H dF_{ij}$) may be severely biased.


\section{Paradoxical Results in the One-Way Layout} \label{owl}

To demonstrate some paradoxical results of rank tests for $d \ge 3$ samples in 
the one-way layout, we consider the vector $\vp = \vW' \vn$ in (\ref{vpdef}) 
of the \np\ effects $p_i$, which are all equal to their mean $\olp_\cdot = 
\frac1d \sumi d p_i$ iff $\sumi d (p_i - \olp_\cdot)^2 = 0$ or in matrix 
notation $\vp' \vT_d \vp = 0$. Here, $\vT_d = \vI_d- \frac1d \vJ_d$ denotes the 
centering matrix, $\vI_d$ the $d$-dimensional unit matrix, and 
$\vJ_d = \veins_d \veins_d'$ the $d\times d$-dimensional matrix of $1$s. Let 
$\vwhp = \int \hH d \vwhF$ denote the plug-in estimator of $\vwhp$ defined in 
(\ref{vechp}). 
In order to detect whether the the $p_i$ are different, we study the \asy\ 
\db\ of $\sqrt{N} \vT_d \vwhp$. This is obtained from the asymptotic 
equivalence theorem (see, e.g., Akritas et al., 1997; Brunner and Puri, 2001, 
2002 or Brunner et al., 2017),
\bqan
 \sqrt{N} \vT_d \vwhp & \appeq\ & \sqrt{N} \vT_d \left[ \volY_\cdot + 
                                  \volZ_\cdot - 2\vp \right] + \sqrt{N} \vT_d 
																	\vp,                  \label{aeqt}
\eqan
where the symbol $\appeq$ denotes \asy\ equivalence. Here, $\volY_\cdot = 
\int H d \vwhF$ and $\volZ_\cdot = \int \hH d \vF$ are vectors of means of 
\ind\ random vectors with expectation $E(\volY_\cdot) = E(\volZ_\cdot) = \vp$. 
It follows from the central limit theorem that $\sqrt{N} \vT_d \left[ 
\volY_\cdot + \volZ_\cdot - 2\vp \right]$ has, \asy ally, a multivariate normal 
\db\ with mean $\vnull$ and \cm\ $\vT_d \vSigma_N \vT_d$, where $\vSigma_N = 
\Cov \left( \sqrt{N} [ \volY_\cdot + \volZ_\cdot]  \right)$ has a quite 
involved structure (for details see Brunner et al., 2017). Obviously, the  
multivariate \db\ is shifted by $\sqrt{N} \vT_d \vp$ from the origin $\vnull$. 
Therefore, we call $\vT_d \vp$ the ``multivariate non-centrality'', 
and a ``univariate non-centrality'' may be quantified by the quadratic form 
$c_p = \vp' \vT_d \vp$. In particular, we have $c_p = 0$ iff $\vT_d \vp = 
\vnull$. The actual (multivariate) shift of the distribution, depending on the 
total sample size $N$, is $\sqrt{N} \vT_d \vp$, and the corresponding 
univariate non-centrality (depending on $N$) is then given by $N \cdot c_p$. 
From these considerations, it should become clear that {\color{black} $N \cdot 
c_p \to \infty$ as $N \to \infty$ if $\vT_d \vp \neq \vnull$.} This defines the 
consistency region of a test based on $\sqrt{N} \vT_d \vwhp$. 

Below, we will demonstrate that for the same vector of \dbs\ $\vF$, the 
non-centrality $c_p = \vp' \vT_d \vp$ may be $0$ in case of equal sample sizes, 
while $c_p$ may be unequal to $0$ in case of unequal sample sizes. Under 
$H_0^F: \vT_d \vF = \vnull$, tests based on $\vwhp$ (such as the Kruskal-Wallis 
test) reject the \yp\ $H_0^F$ with approximately the pre-assigned type-I error 
probability $\alpha$. If, however, the strong hypothesis $H_0^F$ is not true 
then the non-centrality $c_p = \vp' \vT_d \vp$ may be $0$ or unequal to $0$ for 
the same set of \dbs\ $F_1, \ldots, F_d$, since $c_p$ depends on the relative 
samples sizes $n_1/N, \ldots, n_d/N$. This means that for the same set of \dbs\ 
$F_1, \ldots, F_d$ and {\it unequal} sample sizes the $p$-value of the test may 
be arbitrary small if $N$ is large enough. However, the $p$-value for the same 
test may be quite large for the same total sample size $N$ in case of 
{\it equal} sample sizes. Some well-known tests which have this paradoxical 
property are, for example, the Kruskal-Wallis test (1952), the 
Hettmansperger-Norton trend test (1987), and the Akritas-Arnold-Brunner test 
(1997).

As an example, consider the case of $d=3$ \dbs\ where straightforward 
calculations show that 
\ben
 \item in case of equal sample sizes, 
  \bqan
	 p_1 = p_2 = p_3 & \liff\ & w_{21} = w_{32} = 1-w_{31} = w , \label{wijdiff} 
	\eqan
 \item in general, however, 
	\bqan
	 p_1 = p_2 = p_3 & \liff\ & w_{21} = w_{32} = w_{31} = \tfrac12 . 
	 \label{wijid} 
	\eqan
\een

This means that $c_p = \vp' \vT_d \vp = 0$ in case of equal sample sizes, but
$c_p \neq 0$ in case of unequal sample sizes if $w_{21} = w_{32} = 1-w_{31} = w 
\neq \frac12$. 

We note that (\ref{wijid}) follows under the strict \yp\ $H_0^F: F_1 = F_2 = 
F_3$. However, this null hypothesis is not a necessary condition for 
(\ref{wijid}) to hold. For example, if $F_i, \ i=1,2,3$ are symmetric \dbs\ 
with the same center of symmetry then $w_{ri} = \int F_r dF_i = \frac12$ for 
$i=1,2,3$. Thus, in this case, $c_p=0$ is also true for all samples sizes.

An example of discrete \dbs\ generating the \np\ effects $w_{ri}$ in 
(\ref{wijdiff}) is given by the probability mass functions
\bit
 \item $f_1(x) = \frac16$ \ if $x \in \{9,16,17,20,21,22 \}$ and $f_1(x) = 0$ 
       otherwise,
 \item $f_2(x) = \frac16$ \ if $x \in \{13,14,15,18,19,26 \}$ and $f_2(x) = 0$ 
       otherwise,
 \item $f_3(x) = \frac16$ \ if $x \in \{10,11,12,23,24,25 \}$ and $f_3(x) = 0$ 
       otherwise,
\eit
which are derived from some tricky dice (see, e.g., Peterson, 2002). For the 
\dfs\ $F_i(x)$ defined by $f_i(x)$, $i=1,2,3$ above, it is easily seen that 
\bqan
 w_{21} &=& P(X_2 < X_1) = \int F_2 dF_1 = \nfrac7{12} \label{exowlw21} \\
 w_{13} &=& P(X_1 < X_3) = \int F_1 dF_3 = \nfrac7{12} \label{exowlw13} \\
 w_{32} &=& P(X_3 < X_2) = \int F_3 dF_2 = \nfrac7{12} . \label{exowlw32}
\eqan

Thus, $w_{21} = w_{13} = 1 - w_{31} = w_{32} = w$ and the vector of the 
weighted relative effects is given by
\bqan
 \vp &=& \left( \bay{l} p_1 \\ p_2 \\ p_3 \eay \right) \ = \ \vW' \vn \ = \ 
         \frac1N \ \left( \bay{l} 
                            \tfrac12 n_1 + n_3 + (n_2-n_3) w \\[0.5ex]
					 								  n_1 + \tfrac12 n_2 + (n_3-n_1) w \\[0.5ex]
														n_2 + \tfrac12 n_3 + (n_1-n_2) w 
													\eay \right).   \label{exowl3}
\eqan

The weighted relative effects $p_i$ and the resulting non-centralities $c_p$ 
are listed in Table~\ref{cptricky} for equal and some different unequal sample 
sizes.
{\fts
\btab
\caption{\fts Ratios of relative sample sizes $\nfrac{n_i}N$, weighted relative 
effects $p_i$, and the non-centralities for the example of the tricky dice 
where $w=\nfrac{7}{12}$ and the \dbs\ $F_1$, $F_2$, and $F_3$ are fixed.} 
\label{cptricky} \text{ } \\[1ex]

\btb{cccccclllcll} \hline
& & & & & & & \\[-1.5ex] 
Setting & \hspace*{3ex} & $\nfrac{n_1}N$ & $\nfrac{n_2}N$ & $\nfrac{n_3}N$ & 
\hspace*{5ex} & $p_1$ & $p_2$ & $p_3$ & \hspace*{5ex} & $\olp_\cdot$ & $c_p$ 
\\[0.5ex] \cline{1-5} \cline{7-9} \cline{11-12}
& & & & & & & \\[-1.5ex] 
(A) & & $\nfrac13$ & $\nfrac13$ & $\nfrac13$ & & 0.5 & 0.5 & 0.5 & & 0.5 & 0 \\
(B) & & $\nfrac23$ & $\nfrac1{12}$ & $\nfrac14$ & & 0.4861 & 0.4653 & 0.5486 & & 
0.5 & 0.00376 \\
(C) & & $\nfrac14$ & $\nfrac23$ & $\nfrac1{12}$ & & 0.5486 & 0.4861 & 0.4653 & & 
0.5 & 0.00376 \\[0.5ex] \hline
\etb
\etab
}

Since for unequal sample sizes one obtains $c_p \neq 0$, it is only a question 
of choosing the total sample size $N$ large enough to reject the \yp\ $H_0^F: 
F_1 = F_2 = F_3$ by the Kruskal-Wallis test with a probability arbitrary close 
to $1$ while in case of equal sample sizes for $N \to \infty$ the probability 
of rejecting the \yp\ remains constant equal to $\alpha^*$ (close to $\alpha$) 
since in this case, $c_p=0$. It may be noted that in general $\alpha^* \neq 
\alpha$ since the variance estimator of the Kruskal-Wallis \stat\ is computed 
under the strong \yp\ $H_0^F: F_1= F_2 = F_3$, which is obviously not true 
here. Thus, the scaling is not correct, and the Kruskal-Wallis test has a 
slightly different type-I error $\alpha^*$.

For the Hettmansperger-Norton trend test, the situation gets worse since for 
different ratios of sample sizes the \np\ effects $p_1$, $p_2$, and $p_3$ may 
change their order. In setting (B) in Table~\ref{cptricky} we have $p_2 < p_1 < 
p_3$, while in setting (C) we have $p_3 < p_2 < p_1$. Now consider the 
non-centrality of the Hettmansperger-Norton trend test which is a linear rank 
test. Let $\vc = (c_1, \ldots, c_d)'$ denote a vector reflecting the 
conjectured pattern. Then it follows from (\ref{aeqt}) that
\bqan
 \sqrt{N} \vc' \vT_d \vwhp & \appeq\ & \sqrt{N} \vc' \vT_d \left[ \volY_\cdot + 
               \volZ_\cdot - 2\vp \right] + \sqrt{N} \vc' \vT_d	\vp,  
							 \label{nchn}
\eqan
where $c_{HN} = \vc' \vT_d	\vp$ is a univariate non-centrality. If $\vT_d \vF 
= \vnull$ then it follows that $\vT_d	\vp = \vnull$ and $c_{HN} = \vc' \vT_d 
\vp =0$. If, however, $\vT_d \vF \neq \vnull$ then $c_{HN} < 0$ indicates a 
decreasing trend and $c_{HN} > 0$ an increasing trend. In the above discussed 
example, we obtain for setting (B) and for a conjectured pattern of $\vc = 
(1,2,3)'$ for an increasing trend the non-centrality  $c_{HN} = \sumi 3 c_i(p_i 
- \frac12) = \nfrac1{16} > 0$, indeed indicating an increasing trend. For 
setting (C) however, we obtain $c_{HN}=-\nfrac1{12}$, indicating a decreasing 
trend. In case of setting (A) (equal sample sizes), $c_{HN}=0$ since $p_1 = p_2 
= p_3 = \frac12$, and thus indicating no trend. Again it is a question of the 
total sample size $N$ to obtain the decision of a significantly decreasing 
trend for the first setting (B) of unequal sample sizes and for the second 
setting (C) the decision of a significantly increasing trend with a probability 
arbitrary close to 1 for the same \dbs\ $F_1, F_2$, and $F_3$. In the fist 
case, $\sqrt{N} \cdot c_{HN} \to \infty$ for $N \to \infty$ and in the second 
case, $\sqrt{N} \cdot c_{HN} \to - \infty$ for $N \to \infty$. In case of 
equal sample sizes, the \yp\ of no trend is only rejected with a type-I error 
probability $\alpha^{**}$. Regarding $\alpha^{**} \neq \alpha$, a similar 
remark applies as above for the Kruskal-Wallis test.

\section{Paradoxical Results in the Two-Way Layout} \label{twl}

In the previous section, paradoxical decisions by rank tests in case of unequal 
sample sizes were demonstrated for the one-way layout using large sample sizes 
and particular \dbs\ leading to non-transitive decisions. In this section, we 
will show that in two-way layouts paradoxical results are already possible with 
rather small sample sizes and even in simple homoscedastic normal shift models. 
To this end, we consider the simple $2\times 2$-design with two crossed factors 
$A$ and $B$, each with two levels $i=1,2$ for $A$ and $j=1,2$ for $B$. The 
\obss\ $X_{ijk} \sim F_{ij}$, $\kg n_{ij}$, are assumed to be \ind.

The \yps\ of no \np\ effects in terms of the \dfs\ $F_{ij}(x)$ are expressed as
(see Akritas et al., 1997)
\ben
 \item[(1)] no main effect of factor $A$ - $H_0^F(A): F_{11} + F_{12} - F_{21} 
            - F_{22} = 0$
 \item[(2)] no main effect of factor $B$ - $H_0^F(B): F_{11} - F_{12} + F_{21} 
            - F_{22} = 0$
 \item[(3)] no interaction $AB$ - $H_0^F(AB): F_{11} - F_{12} - F_{21} + F_{22} 
            = 0$,
\een
where in all three cases, $0$ denotes a function which is identical $0$.

Let $\vF = (F_{11}, F_{12}, F_{21}, F_{22})'$ denote the vector of the \dfs.
Then the three \yps\ formulated above can be written in matrix notation as 
$H_0^F(\vc): \vc' \vF = 0$, where $\vc = \vc_A = (1,1,-1,-1)'$ generates the 
hypothesis for the main effect $A$, $\vc = \vc_B = (1,-1,1,-1)'$ for the main 
effect $B$, and $\vc = \vc_{AB} = (1,-1,-1,1)'$ for the interaction $AB$.

For testing these \yps, Akritas et al. (1997) derived rank procedures based on 
the \stat\ 
\bqan
 T_N(\vc) &=& \sqrt{N} \vc' \vwhp \ = \ \frac1{\sqrt{N}} \ \vc' \volR_\cdot,  
\label{tnvc}
\eqan
where $\volR_\cdot = (\olR_{11\cdot}, \olR_{12\cdot}, \olR_{21\cdot}, 
\olR_{22\cdot})'$ denotes the vector of the rank means $\olR_{ij \cdot}$ within 
the four samples. They showed that under the \yp\ $H_0^F(\vc)$, the \stat\ 
$T_N(\vc)$ has, \asy ally, a normal \db\ with mean $0$ and variance
\bqan
 \sigma_0^2 &=& \sumi 2 \sumj 2 \frac{N}{n_{ij}} \sigma_{ij}^2 , \label{sigijq} 
\eqan
where the unknown variances $\sigma_{ij}^2$ (see Akritas et al., 1997, for 
their explicit form) can be consistently estimated 
by 
\bqan
\frac1{N^2} S_{ij}^2 &=& \frac1{N^2(n_{ij}-1)} \sumk {n_{ij}} (R_{ijk} - 
                         \olR_{ij \cdot})^2   \label{sijq}
\eqan												
and $\hsigma_0^2 = N \ \sumi 2 \sumj 2 S_{ij}^2 / n_{ij}$. For small sample 
sizes, the null \db\ of $L_N(\vc) = T_N(\vc) / \hsigma_0$ can be approximated 
by a $t_f$-\db\ with estimated degrees of freedom 
\bqan
 \hf &=& \frac{S_0^4}{\sumi 2 \sumj 2 (S_{ij}^2/n_{ij})^2 / (n_{ij}-1)} ,
         \label{fhat}
\eqan
where $S_0^2 = \sumi 2 \sumj 2 S_{ij}^2/n_{ij}$. The non-centrality of 
$T_N(\vc)$ is given by $c_T = \vc' \vp$, and under the restrictive null \yp\ 
$H_0^F(\vc): \vc' \vF = 0$ it follows that $\vc_T' \vp = 0$.


To demonstrate a paradoxical result, we assume that the \obss\ $X_{ijk}$ are 
coming from the normal \dbs\ $N(\mu_{ij}, \tau^2)$ with equal standard 
deviations $\tau=0.4$ and expectations $\vmu = (\mu_{11}, \mu_{12}, \mu_{21}, 
\mu_{22})' = (10,9,9,8)'$. From the viewpoint of linear models, there is a main 
effect $A$ of $\vc_A' \vmu = \mu_{11} + \mu_{12} - \mu_{21} - \mu_{22} = 2$, a 
main effect $B$ of $\vc_B' \vmu = \mu_{11} - \mu_{12} + \mu_{21} - \mu_{22} = 
2$, and no $AB$-interaction since $\vc_{AB}' \vmu = \mu_{11} - \mu_{12} - 
\mu_{21} + \mu_{22} = 0$. Since this is a homoscedastic linear model, the 
classical ANOVA should reject the \yps\ $H_0^\mu(\vc_A): \vc_A' \vmu = 0$ and 
$H_0^\mu(\vc_B): \vc_B' \vmu = 0$ with a high probability if the total sample 
size is large enough. In contrast to that, the \yp\ $H_0^\mu(\vc_{AB}): 
\vc_{AB}' \vmu = 0$ of no interaction is only rejected with the pre-selected 
type-I error probability $\alpha$. The non-centralities are given by $c_A^\mu = 
\vc_A' \vmu = 2$, $c_B^\mu = \vc_B' \vmu = 2$, and $c_{AB}^\mu = \vc_{AB}' \vmu 
= 0$. The following two settings of samples sizes $n_{ij}$ are considered:
\ben
 \item[(1)] $n_{11} = 10, \ n_{12} = 20, \ n_{21} = 20, \ n_{22} = 50$, \sep - 
            (unbalanced)
 \item[(2)] $n_{11} = n_{12} = n_{21} = n_{22} = 25$, \sep - (balanced).
\een

First we demonstrate that the empirical characteristics of the two data sets, 
which are sampled from the same \dbs, are nearly identical. Thus, potentially 
different results could not be explained by substantially different empirical \dbs\ 
obtained by an ``unhappy randomization''. 
The results of the comparisons are listed in Table~\ref{2x2.comp}.
\text{ } \\[-2ex]
\btab
{\fts
\caption{\fts Comparison of the empirical distributions in the four factor level 
combinations $A_1 B_1$, $A_2 B_1$, $A_1 B_2$, and $A_2 B_2$ 
sampled from homoscedastic normal \dbs\ with 
$\mu_{11} = 10, \mu_{12} = \mu_{21} = 9, \mu_{22}=8$, 
and standard deviation $\tau=0.4$. 
Within each factor level combination $A_i B_j$, $i,j=1,2$, 
the unadjusted $p$-values of a $t$-test comparing the location of 
balanced and unbalanced samples, 
more specifically testing $\mu_{ij}(\text{balanced}) = \mu_{ij}(\text{unbalanced})$, 
are listed in the last column.} \label{2x2.comp} 
\bcen
\btb{ccccccccccc} \hline 
& & & & &  \\[-1.5ex]
Level & \hspace*{3ex} & \mc{3}{c}{Means} & \hspace*{3ex} & 
\mc{3}{c}{Standard Deviations} & \hspace*{3ex} & $t$-Tests \\[0.3ex]
\cline{3-11} 
& & & & &  \\[-1.7ex]
Combinations & & Balanced & \hspace*{1ex} & Unbalanced & & Balanced & 
\hspace*{1ex} & Unbalanced & & $p$-Values \\[0.4ex]
\cline{1-1} \cline{3-5} \cline{7-9} \cline{11-11} 
& & & & &  \\[-1.5ex]
$A_1 B_1$ & & 10.07& & 9.91 & & 0.311 & & 0.314 & & 0.178 \\
$A_2 B_1$ & & 9.04 & & 8.95 & & 0.380 & & 0.313 & & 0.409 \\
$A_1 B_2$ & & 9.05 & & 8.99 & & 0.408 & & 0.480 & & 0.684 \\
$A_2 B_2$ & & 8.07 & & 8.02 & & 0.371 & & 0.359 & & 0.562 \\[0.4ex] \hline
\etb
\ecen
} 
\etab

We apply the classical ANOVA $F$-\stat\ and the rank statistic $L_N(\vc) = 
T_N(\vc) / \hsigma_0$ to the same simulated data sets from Table~\ref{2x2.comp} 
and compare the results for testing the three \yps\ $H_0^F(A)$, $H_0^F(B)$, and 
$H_0^F(AB)$. We note that $H_0^F(\vc): \vc' \vF = 0 \RA\ H_0^\mu(\vc): \vc' 
\vmu = 0$ and $\vc' \vmu \neq 0 \RA\ \vc' \vF \neq 0$. The decisions for 
$H_0^F(A)$, $H_0^F(B)$, and  $H_0^F(AB)$ obtained by the ANOVA as well as by 
the rank tests based on $L_N(\vc)$ are identical in all cases in the balanced 
setting. In the unbalanced setting, all decisions obtained by the parametric 
ANOVA are comparable to those in the balanced setting. However, the decision on 
the interaction $AB$ based on the rank test is totally different from that 
obtained by the parametric ANOVA, as well as that obtained for the rank test in 
the balanced setting. The results are summarized in Table~\ref{tab2x2}.

\btab[t]
{\fts
\caption{\fts Comparison of the results obtained by an ANOVA and by the rank 
test $L_N(\vc)$ in case of a balanced (left) and an unbalanced (right) 
$2\times 2$-design. Surprising is the fact that in the balanced case, the 
decisions of both procedures coincide while in the unbalanced case 
the decisions for testing the interaction $AB$ are totally opposite.} 
\label{tab2x2} 
\bcen
\btb{ccrccrccrccrc} \hline
& & & & & \\[-1.5ex]
& \hspace*{3ex} & \mc{5}{c}{Balanced} & \hspace*{3ex} & \mc{5}{c}{Unbalanced}\\ 
& & \mc{5}{c}{$n_{11}=n_{12}=n_{21}=n_{22}=25$} & & 
\mc{5}{c}{$n_{11}=10, n_{12}=20, n_{21}=20, n_{22} = 50$} 
\\[0.5ex] \cline{3-13}
& & & & & \\[-1.5ex]
& & \mc{2}{c}{ANOVA} & & \mc{2}{c}{Rank Test} & \hspace*{3ex} & 
    \mc{2}{c}{ANOVA} & & \mc{2}{c}{Rank Test}\\[0.5ex]\cline{3-7}\cline{9-13}
& & & & & \\[-1.5ex]
Effect & & $F$ & $p$-Value & \hspace*{1ex} & $L_N^2(\vc)$ & $p$-Value & & $F$ & 
$p$-Value & \hspace*{1ex} & $L_N^2(\vc)$ & $p$-Value 
\\[0.5ex] \cline{1-1} \cline{3-4} \cline{6-7} \cline{9-10} \cline{12-13}
& & & & & \\[-1.5ex]
$A$ & & 184.43 & $<10^{-4}$ & & 168.19 & $<10^{-4}$ & & 120.99 & $<10^{-4}$ & & 
150.07 & $<10^{-4}$ \\
$B$ & & 182.60 & $<10^{-4}$ & & 169.32 & $<10^{-4}$ & & 111.48 & $<10^{-4}$ & & 
143.55 & $<10^{-4}$ \\
$AB$ & & 0.12 & 0.7317 & & 0.00 & 0.9541 & & 0.01 & 0.9112 & & 7.99 & 0.0065 
\\[0.5ex] \hline
\etb
\ecen
}
\etab

On the surface, the difference of the decisions in the unbalanced case could be 
explained by the fact that the \np\ \yp\ $H_0^F(AB)$ and the parametric \yp\
$H_0^\mu(AB)$ are not identical and that this particular configuration of 
normal \dbs\ falls into the gap between $H_0^F(AB)$ and $H_0^\mu(AB)$. That is, 
here $H_0^\mu(AB)$ is true, but $H_0^F(AB)$ is not. It is surprising, however, 
that this explanation does not hold for the balanced case. The difference of 
the two $p$-values $0.9541$ and $0.0065$ in Table~\ref{tab2x2} calls for an 
explanation.

The reason becomes clear when computing the vector $\vp$ in (\ref{vpdef}) for 
this particular example of the $2\times2$-design. To avoid fourfold indices we 
re-label the \dbs\ $F_{11}, F_{12}, F_{21}$, and $F_{22}$ as $F_1, F_2, F_3$, 
and $F_4$, respectively, and the sample sizes accordingly as $n_1, n_2, n_3$, 
and $n_4$. In the example, $F_1 = N(10, \tau^2)$, $F_2 = F_3 = N(9, \tau^2)$, 
and $F_4 = N(8, \tau^2)$, where $\tau = 0.4$. Thus, the probabilities $w_{ri} = 
\int F_r d F_i$ of the pairwise comparisons are
\bqa
 w = w_{12} = w_{13} = w_{24} = w_{34} &=& \Phi \left(
                                         \frac{-1}{\tau \sqrt{2}} \right) \ = \ 
																				 0.0392, \\
 w_{23} = w_{32} &=& \frac12, \\ 
 v = w_{14} &=& \Phi \left(\frac{-\sqrt{2}}{\tau} \right) \ \approx\ \ 0. 
\eqa

Finally, by observing $w_{ri} = 1 - w_{ir}$, we obtain
\bqa
 \vp &=& \vW' \vn \ = \ \frac1N 
                      \left( \bay{l} 
											 \frac12 n_1 + n_2 + n_3 + n_4 - (n_2+n_3) w - n_4 v 
											  \\[1ex]
											 \frac12 (n_2+n_3) + n_4 + (n_1-n_4) w \\[1ex]
											 \frac12 (n_2+n_3) + n_4 + (n_1-n_4) w \\[1ex]
											 \frac12 n_4 + (n_2+n_3) w + n_1 v
											 \eay
											\right) ,
\eqa
and the nonparametric $AB$-interaction is described by
\bqan
c_{AB}^p \ = \ \vc_{AB}' \ \vp &=& p_1 - p_2 - p_3 + p_4  \label{cabp} \\
				&=& \frac{n_1-n_4}N \left( \tfrac12 - 2 w + v \right) . \nnr\ 
\eqan

In this example, we obtain for equal samples sizes $c_{AB}^p = 0$, while for 
$\vn = (10, 20, 20, 50)$ we obtain $c_{AB}^p= - \frac25 \left(\frac12 - 2w + v 
\right) \approx - 0.1686$ and $\sqrt{N} c_{AB}^p \approx - 1.686$. This 
explains the small $p$-value for unequal sample sizes in the example.

\section{Explanation of the Paradoxical Results} \label{expl}
\subsection{Unweighted Effects and Pseudo-Ranks} \label{pseudo}

The simple reason for the paradoxical results is the fact that even when all 
distribution functions underlying the observations are specified, the 
consistency regions $\vC' \vp \neq \vnull$ of the rank tests based on $\vC' 
\vwhp$ are not fixed. Indeed, the consistency regions are defined by the 
weighted nonparametric relative effects $p_i$ which are not fixed model 
quantities by which \yps\ could be formulated or for which \cis\ could be 
reasonably computed since the $p_i$ themselves generally depend on the sample 
sizes $n_i$. 

Thus, it appears reasonable to define different \np\ effects which are fixed 
model quantities not depending on sample sizes. To this end let $G = \frac1d 
\sumr d F_r$ denote the unweighted mean distribution, and let $\psi_i = \int G d F_i$. 
Easily interpreted, this \np\ effect $\psi_i$ measures an effect of the \db\ 
$F_i$ with respect to the unweighted mean \db\ $G$ and is therefore a 
``fixed relative effect''. As $\psi_i = \frac1d \sumr d w_{ri}$ is the mean of the 
pairwise \np\ effects $w_{1i}, \ldots, w_{di}$, it can be written in vector 
notation as the vector of row means of $\vW'$, that is,
\bqan
 \vpsi &=& \int G d \vF \ = \ \vW' \cdot\ \tfrac1d \veins_d \ = \ \left( 
           \bay{ccc} 
           w_{11}, & \cdots & w_{d1} \\ \vdots & \ddots & \vdots \\ w_{1d}, & 
					 \cdots & w_{dd} 
					 \eay 
					 \right) \cdot \frac1d \veins_d \ = \ \left( \bay{c} \psi_1 \\ \vdots 
					 \\ \psi_d \eay \right). \label{vpsidef}
\eqan

The fixed relative effects $\psi_i$ can be estimated consistently by the simple 
plug-in estimator
\bqan 
 \hpsi_i &=& \int \hG d \hF_i \ = \ \frac1N \left(\olR_{i \cdot}^\psi - 
             \tfrac12 \right),  \label{psiest}
\eqan
where $\hG = \frac1d \sumr d \hF_r$ denotes the unweighted mean of the 
empirical \dbs\ $\hF_1, \ldots, \hF_d$, and $\olR_{i \cdot}^\psi = \frac1{n_i} 
\sumk {n_i} R_{ik}^\psi$ the mean of the so-called {\it pseudo-ranks}
\bqan
 \text{ps-rank}(X_{ik}) \ = \ R_{ik}^\psi &=& \frac12 + N \hG(X_{ik}) \ = \ 
 \frac12 + \frac{N}d \sumr d \frac1{n_r} \sumell {n_r} c(X_{ik} - X_{r \ell}). 
 \label{rikpsidef}
\eqan

Finally, the estimators $\hpsi_i$ are arranged in the vector 
\bqan
 \vwhpsi &=& \left( \bay{c} \hpsi_1 \\ \vdots \\ \hpsi_d \eay \right) \ = \
             \int \hG d \vwhF \ = \ \frac1N \left( \volR_\cdot^\psi - \tfrac12 
             \veins_d \right),                        \label{vechpsi}
\eqan
where $\volR_\cdot^\psi = (\olR_{1 \cdot}^\psi, \ldots, \olR_{d \cdot}^\psi)'$ 
is the vector of the pseudo-rank means $\olR_{i \cdot}^\psi$.

It may be noted that the pseudo-ranks $R_{ik}^\psi$ have similar properties as 
the ranks $R_{ik}$. The properties given below follow from the definitions of 
the ranks and pseudo-ranks and by some straightforward algebra. 
\begin{lemma} \label{rpsrcomp}
Let $X_{ik}$ denote $N = \sumi d n_i$ \obss\ arranged in $\ig d$ groups each 
involving $n_i$ \obss. Then, for $i,\rg d$ and $\kg n_i, \ellg n_r$,
\ben
 \item If $X_{ik} < X_{r\ell}$ then $R_{ik} < R_{r\ell}$ and $R_{ik}^\psi < 
       R_{r\ell}^\psi$.
 \item If $X_{ik} = X_{r\ell}$ then $R_{ik} = R_{r\ell}$ and $R_{ik}^\psi = 
       R_{r\ell}^\psi$. 			
 \item $\frac1N \sumi d \sumk {n_i} R_{ik} \ = \ \frac1d \sumi d \frac1{n_i} 
       \sumk {n_i} R_{ik}^\psi \ = \ \frac{N+1}2$.			
 \item If $m(u)$ is a strictly monotone transformation of $u$ then,
   \ben
		\item[(a)] $R_{ik} = \text{rank}(X_{ik}) = \text{rank}(m(X_{ik}))$
		\item[(b)] $\dl R_{ik}^\psi = \text{ps-rank}(X_{ik}) = \text{ps-rank} 
		           \left(m(X_{ik}) \right)$.
	 \een
 \item $1 \le R_{ik} \le N$.
 \item $\dl \frac{d+1}{2d} \le \ \frac12 + \frac{N}{2dn_i} \ \le \ R_{ik}^\psi 
       \ \le N +  \frac12 - \frac{N}{2dn_i} \ \le N + \frac{d-1}{2d}$
 \item If $n_i \equiv n, \ \ig d$, then $R_{ik} = R_{ik}^\psi$
 \item Let $X_{ik}$, $\ig d; \kg {n_i}$ be independent and identically 
       distributed random variables, then 
			 \bqa 
			   E \left( R_{ik} \right) \ = \ E \left( R_{ik}^\psi  \right) &=& 
				    \frac{N+1}2 \ .
			 \eqa			
\een
\end{lemma}

\subsection{Consistency Regions of Pseudo-Rank Procedures} \label{pseudo.cons}

As a solution to the paradoxical results discussed in Sections~\ref{owl} and 
\ref{twl}, we demonstrate that replacing the ranks $R_{ik}$ with the pseudo-ranks 
$R_{ik}^\psi$ leads to procedures that do not have these undesirable 
properties. 
The main reason is that pseudo-rank procedures are based on the 
(unweighted) relative effects $\psi_i$ which are fixed model quantities 
by which hypotheses can be formulated and for which \cis\ can be derived. 
In case of equal sample sizes $n_i \equiv n$, $\ig d$, 
we do not obtain paradoxical results since in this case 
ranks and pseudo-ranks coincide,
$R_{ik} = R_{ik}^\psi$ (see Lemma~\ref{rpsrcomp}). 

Pseudo-rank based inference procedures are obtained in much the same way as the 
common rank procedures, by using relations (\ref{psiest}) and (\ref{rikpsidef}),
which generally means substituting ranks by the corresponding pseudo-ranks. 
Only in the computation of the \cis\ (see Section~\ref{expl.confi}), there is a 
minor change in the variance estimator. For details we refer to Brunner, 
Bathke, and Konietschke (2018, Result~4.16 in Section~4.6.1). 

Below, we examine the behavior of pseudo-rank procedures in the
situations where the use of rank tests led to paradoxical results. 

\ben
 \item For testing the hypothesis $H_0^F: F_1 = F_2 = F_3 \ \liff\ \ \vT_3 \vF 
       = \vnull$ in case of the tricky dice in the example in Section~\ref{owl}, 
			 one obtains for the non-centrality of the Kruskal-Wallis statistic, 
			 when the ranks $R_{ik}$ are replaced with the pseudo-ranks 
			 $R_{ik}^\psi$, the value $c_\psi = \vpsi' \ \vT_3 \ \vpsi = 0$. This is 
			 easily seen from (\ref{vpsidef}), since in this case it follows from 
			 (\ref{exowlw21}), (\ref{exowlw13}), and (\ref{exowlw32}) that $\vpsi = 
			 \vW' \frac13 \veins_3 = \frac12 \veins_3$ and $c_\psi = \vpsi' \ \vT_3 \ 
			 \vpsi = \frac12 \veins_3' \ \vT_3 \ \frac12 \veins_3 = 0$ by noting that 
			 $\vT_3$ is a contrast matrix and thus, $\vT_3 \veins_3 = \vnull$. 
 \item The non-centrality for the Hettmansperger-Norton trend test when substituting 
       ranks with pseudo-ranks becomes $c_{HN}^\psi = \vc'\ \vT_3 \ \vpsi 
			 = \vc' \ \vT_3 \ \vW' \ \frac13 \veins_3 = 0$ for all trend alternatives 
			 $\vc = (c_1, c_2, c_3)'$.
 \item In the two-way layout, we reconsider the example of the four shifted 
       normal \dbs\	we obtain the unweighted relative effects $\psi_{ij}$ 
			 \bqa
			  \vpsi \ = \ \vW' \cdot\ \tfrac14 \veins_4 &=& \frac14 \left( 
				  \bay{l} 
					 \nfrac72 - 2w - v \\
					 2 \\ 2 \\ \nfrac12 + 2w + v 
					\eay 
					\right) 
			 \eqa
and the non-centrality of the statistic for the interaction 
\bqan
 c_{AB}^\psi &=& \psi_1 - \psi_2 - \psi_3 + \psi_4 \ = \ 0. \label{cabpsi}
\eqan 

In the analysis of the unbalanced data in Table~\ref{tab2x2} using pseudo-ranks, 
one obtains for the interaction the statistic $L_N^2(\vc_{AB}) = 0.08$ and the 
resulting $p$-value is $0.7832$ which agrees with the result for the balanced 
sample from the same distributions in Table~\ref{tab2x2}.
\een

In all cases, paradoxical results obtained by changing the ratios of the sample 
sizes cannot occur since the non-centralities $c_\psi$, $c_{HN}^\psi$, and 
$c_{AB}^\psi$ are equal to $0$ for all constellations of the relative sample 
sizes.

In case of $d=2$ samples, it is easily seen that $p_2 - p_1 = \psi_2 - \psi_1 = 
p = \int F_1 d F_2$ which does not depend on sample sizes. Thus, paradoxical 
results for rank-based tests -- 
such as presented in the previous sections -- 
can only occur for $d \ge 3$ samples.

\subsection{Confidence Intervals} \label{expl.confi}

Here we briefly explain the details on how to compute confidence intervals for 
fixed nonparametric effects which have an intuitive and easy to understand 
interpretation. This shall be demonstrated by means of the example involving 
the four shifted normal \dbs\ considered in Section~\ref{twl}. 

The quantity $\hpsi_i$ estimates the probability that a randomly drawn 
observation from the mean distribution $G = \frac14 \sumi 4 F_i$ is smaller 
than a randomly drawn observation from distribution $F_i$ (plus $\frac12$ times 
the probability that they are equal). The estimator is obtained from (\ref{psiest}), 
and the limits of the confidence interval are obtained from formula (25) in 
Brunner et al. (2017) or Result~4.16 in Chapter~4 of Brunner, Bathke, and 
Konietschke (2018). The results are listed in Table~\ref{expl.tab}. \\[-3ex]

\btab[h]
{\fts
\caption{\fts Estimates and two-sided $95\%$-\cis\ $[\psi_{i,L}, \psi_{i,U}]$ 
for the (unweighted) nonparametric relative effects $\psi_i$ of the two data sets 
considered in Table~\ref{tab2x2}. 
The results for the balanced case (upper part) and the 
unbalanced case (lower part) are quite similar.} \label{expl.tab} 
\text{ } \\[-7ex]
	\bcen
	\btb{ccccccc} \\ \hline
	& & & & &  \\[-1.5ex]
	& & \mc{5}{c}{Equal Sample Sizes $n_1=n_2=n_3=n_4=25$} \\[0.4ex] \hline 
	& & & & &  \\[-1.5ex]
	\mc{1}{c}{Sample} & \hspace*{3ex} & \mc{1}{c}{Lower Limit $\psi_{i,L}$} & 
	\hspace*{1ex} & \mc{1}{c}{Estimator $\psi_i$}  & \hspace*{1ex} & 
	\mc{1}{c}{Upper Limit $\psi_{i,U}$} \\[0.4ex]
	\cline{1-1} \cline{3-7} 
	& & & & & \\[-1.5ex]
	1 & & 0.84 & & 0.86 & & 0.88 \\
	2 & & 0.45 & & 0.50 & & 0.55 \\
	3 & & 0.45 & & 0.50 & & 0.55 \\
	4 & & 0.12 & & 0.14 & & 0.16 \\[0.4ex] \hline 
	& & & & & \\[-1.5ex]
	& & \mc{5}{c}{Unequal Sample Sizes $n_1=10, n_2=n_3=20, n_4=50$} 
	\\[0.4ex] \hline
	& & & & & \\[-1.5ex]
	\mc{1}{c}{Sample} & \hspace*{3ex} & \mc{1}{c}{Lower Limit $\psi_{i,L}$} & 
	\hspace*{1ex} & \mc{1}{c}{Estimator $\psi_i$}  & \hspace*{1ex} & 
	\mc{1}{c}{Upper Limit $\psi_{i,U}$} \\[0.4ex]
	\cline{1-1} \cline{3-7} 
	& & & & & \\[-1.5ex]
	1 & & 0.85 & & 0.86 & & 0.88 \\
	2 & & 0.45 & & 0.51 & & 0.57 \\
	3 & & 0.43 & & 0.48 & & 0.53 \\
	4 & & 0.13 & & 0.15 & & 0.16 \\[0.4ex] \hline 
	\etb
	\ecen
	} 
\etab
\text{ } \\[-3ex]

Obviously, the results for equal and unequal sample sizes are nearly identical. 
Let samples 1 and 2 refer to the factor level combinations 
$A_1B_1$ and $A_2B_1$ within level $B_1$ 
and samples 3 and 4 to the factor level combination $A_1B_2$ 
and $A_2B_2$ within level $B_2$, respectively. 
The differences $\hpsi_1 - \hpsi_2$ and 
$\hpsi_3 - \hpsi_4$ between factor levels $A_1$ and $A_2$ of the estimates in 
factor level $B_1$ and in factor level $B_2$ are identical in the balanced case 
($\hpsi_1 - \hpsi_2 = 0.36$, $\hpsi_3 - \hpsi_4 = 0.36$) and nearly identical 
($\hpsi_1 - \hpsi_2 = 0.35$, $\hpsi_3 - \hpsi_4 = 0.33$) in the unbalanced 
case. These results do not indicate any interaction. 
This finding agrees with the analysis of the balanced data in 
Table~\ref{tab2x2} when using ranks. 

For the weighted relative effects $p_i$ in (\ref{statmod.pi}), the estimators 
$\hp_i$ based on ranks, as well as $95\%$-intervals for $p_i$ are listed in 
Table~\ref{tab2x2.pi}. The $95\%$-intervals for $p_i$, however, cannot strictly 
be interpreted as \cis\ since the quantities $p_i$ depend on the relative 
sample sizes, unless the design is balanced. In case of equal sample sizes we 
have $p_i=\psi_i$. Therefore only the results for the unbalanced case are 
listed in Table~\ref{tab2x2.pi}. The $95\%$-intervals for $p_i$ are obtained 
from formulas (1.14) - (1.17) in Brunner and Puri (2001) and formulas (2.2.33) 
- (2.2.35) in Section~2.2.7 in Brunner and Munzel (2013). 

Comparing the differences of the estimates within the two factor levels $B_1$ 
and $B_2$ using {\color{black} the usual} ranks, we obtain $\hp_1 - \hp_2 = 0.25$ 
and $\hp_3 - \hp_4 = 0.41$ which are quite different. Such a difference 
indicates an interaction effect, which agrees with the results of the analysis 
of the unbalanced data in Table~\ref{tab2x2}. Moreover, the corresponding 
$95\%$-intervals for the weighted effects $p_i$ are totally different from the 
$95\%$-\cis\ for the unweighted effects $\psi_i$ in Table~\ref{expl.tab}. They 
even do not overlap with the \cis\ for the $\psi_i$. Thus, only the 
nonparametric effects $\psi_i$ offer the possibility of computing estimates and 
\cis\ for fixed model quantities with an intuitive interpretation. \\[-3ex] 
\btab[h]
{\fts
\caption{\fts Estimates and $95\%$-intervals $[p_{i,L}, p_{i,U}]$ for the 
weighted effects $p_i$ of the data set in the unbalanced case in 
Table~\ref{2x2.comp}. The results for the balanced case are identical to those 
in Table~\ref{expl.tab} since $p_i = \psi_i$ for equal sample sizes.} 
\label{tab2x2.pi}
\text{ } \\[-6ex]
\bcen
	\btb{ccccccc} \\ \hline
	& & & & &  \\[-1.5ex]
  & & \mc{5}{c}{Unequal Sample Sizes $n_1=10, n_2=n_3=20, n_4=50$} 
	\\[0.4ex] \hline
	& & & & & \\[-1.5ex]
	\mc{1}{c}{Sample} & \hspace*{3ex} & \mc{1}{c}{Lower Limit $p_{i,L}$} & 
	\hspace*{1ex} & \mc{1}{c}{Estimate $\hp_i$}  & \hspace*{1ex} & 
	\mc{1}{c}{Upper Limit $p_{i,U}$} \\[0.4ex]
	\cline{1-1} \cline{3-7} 
	& & & & & \\[-1.5ex]
	1 & & 0.93 & & 0.94 & & 0.95 \\
	2 & & 0.62 & & 0.69 & & 0.74 \\
	3 & & 0.63 & & 0.68 & & 0.72 \\
	4 & & 0.25 & & 0.27 & & 0.28 \\[0.4ex] \hline 
	\etb
\ecen
}
\etab

\section{Cautionary Note for Sub-Group Analysis} \label{cauno.sg}

The fact that the (weighted) relative effects in (\ref{vpdef}), which are 
estimated by the ranks, depend on the relative sample sizes $n_i/N$ is 
important in sub-group analyses, for example in clinical trials. Here, 
typically, the total group of patients is quite large while the sub-group may 
be small. The design of such a trial can be considered as a $2\times 2$-design 
where $F_{11}$ and $F_{12}$ are the distributions of the outcome for the 
patients  who received treatment 1 (e.g., standard treatment with distribution 
$F_{11}$) or treatment 2 (experimental treatment with distribution $F_{12}$), 
without the particular sub-group of interest. The corresponding larger sample 
sizes, $n_{11}$ and $n_{12}$ are approximately the same when using an 
appropriate randomization procedure. The smaller sample sizes $n_{21}$ and 
$n_{22}$ in the sub-group may be equal or quite different depending on whether 
or not the randomization was also stratified for the sub-group.

The main question regarding the sub-group in this design is whether the 
treatment effect is the same as in the population of patients without the 
sub-group. Here we consider only the case where the sub-group is 
known in advance and is not identified on the  basis of the data. 
In order to demonstrate the advantage of using the unweighted relative 
effects in (\ref{vpsidef}) estimated by the pseudo-ranks, 
we consider the homoscedastic two-way design from Section~\ref{twl}. 
To this end, let $\vc_{AB} = (1, -1, -1, 1)'$ denote the contrast vector 
for the interaction in the $2 \times 2$-design 
as given in Section~\ref{twl} and let $\vp$ denote the vector of the weighted 
relative effects estimated by $\vwhp$ using the ranks. Moreover, let $\vpsi$ 
denote the vector of the unweighted relative effects estimated by $\vwhpsi$ 
using pseudo-ranks, and finally let $\vmu$ denote the vector of the 
expectations. Then the corresponding non-centralities are given by $c_{AB}^p = 
\vc_{AB}' \vp$, $c_{AB}^\psi = \vc_{AB}' \vpsi$, and $c_{AB}^\mu = \vc_{AB}' 
\vmu$. The total sample size $N$ is increased such that the sample sizes 
$n_{21} = n_{22} = 50$ in the sub-group remain constant while the samples sizes 
$n_{11} = n_{12}$ obtained from the population without the sub-group increases 
from $50$ to $2000$. For the interaction term of the procedure based on ranks, 
we also consider the quantity $\sqrt{N} c_{AB}^p = \sqrt{N} \vc_{AB}' \vp$ 
describing the shift of the \asy\ multivariate normal \db. It is noteworthy 
that with increasing total sample size $N$ not only $\sqrt{N} c_{AB}^p$ is 
increasing but also the basic non-centrality $c_{AB}^p$ increases while the 
corresponding non-centralities $c_{AB}^\psi$ and $c_{AB}^\mu$ of the two other 
procedures remain constant equal to $0$. Note that for equal sample sizes, all 
non-centralities are $0$. The results are listed in Table~\ref{appl.sgct.tab}.
{\fts
\btab[t]
\caption{\fts Non-centralities $c_{AB}^p$, $c_{AB}^\psi$, and $c_{AB}^\mu$ of 
the interaction statistics in the $2\times 2$-design with increasing sample 
sizes in only one stratum. This reflects the situation of a sub-group analysis 
where the sample size in the sub-group is small compared to the total 
population. For the procedure based on ranks, the quantity $\sqrt{N} c_{AB}^p$ 
is also listed to demonstrate the shift of its \asy\ multivariate normal \db\ 
generated only by the unequal sample sizes in the two strata. The 
non-centralities $c_{AB}^\psi$ and $c_{AB}^\mu$ of the two \stats\ based on the 
pseudo-ranks or on the expectations are both equal to $0$, 
independently of the sample sizes. } \label{appl.sgct.tab}
{\fts
\bcen
\btb{rcrcccccccr} \\ \hline
& & & & & & & & & \\[-1.5ex]
\mc{3}{c}{Sample Sizes} & \hspace*{5ex} & \mc{7}{c}{Non-Centralities} \\[0.4ex] 
\cline{1-3} \cline{5-11}
& & & & & & & & & \\[-1.5ex]
$n_{11} = n_{12}$ & & $n_{21} = n_{22}$ & & $c_{AB}^\mu$ & & $c_{AB}^\psi$ & & 
$c_{AB}^p$ & & $\sqrt{N} c_{AB}^p$ \\[0.4ex] 
\cline{1-3} \cline{5-11}
& & & & & & & & & \\[-1.5ex]
  50 \hspace*{2ex} & &  50 \hspace*{2ex} & & 0 & & 0 & & 0 & &  \mc{1}{c}{0} \\
 100 \hspace*{2ex} & &  50 \hspace*{2ex} & & 0 & & 0 & & 0.071 & &  1.22 \\	
 200 \hspace*{2ex} & &  50 \hspace*{2ex} & & 0 & & 0 & & 0.127 & &  2.84 \\	
 300 \hspace*{2ex} & &  50 \hspace*{2ex} & & 0 & & 0 & & 0.151 & &  4.00 \\	
 400 \hspace*{2ex} & &  50 \hspace*{2ex} & & 0 & & 0 & & 0.165 & &  4.94 \\	
 500 \hspace*{2ex} & &  50 \hspace*{2ex} & & 0 & & 0 & & 0.173 & &  5.74 \\	
1000 \hspace*{2ex} & &  50 \hspace*{2ex} & & 0 & & 0 & & 0.191 & &  8.77 \\	
2000 \hspace*{2ex} & &  50 \hspace*{2ex} & & 0 & & 0 & & 0.201 & & 12.89 
\\[0.4ex] \hline	
\etb
\ecen
}
\etab
}

Analyses of this type are very common in biostatistical practice. In 
particular, the sub-group's proportion of the total sample size is often very 
small, as, for example. in case of rare disease studies with large control 
groups. {\color{black} The above example further illustrates the need for 
changing the usage of classical rank procedures to more favorable pseudo-rank 
procedures.}

\section{Discussion and Conclusions} \label{discon}

We have demonstrated that inference based on ranks may lead to paradox results
in the case of at least three samples and an unbalanced design. These results 
may occur in rather natural situations and are not restricted to artificial 
configurations of the population distributions. 

The reason for this behavior of rank-based tests (including some classical and 
frequently used tests) can be explained by the non-centralities of the test 
statistics which are functions of weighted nonparametric relative effects. 
{\color{black} Moreover, it is worth to note that neither pairwise nor stratified 
rankings provide a solution to the problem. In fact, it would only 
cause new problems. For a discussion of these types of rankings see Brunner et 
al., (2018).} As a remedy, the use of unweighted nonparametric relative 
effects, and thus, the use of tests based on pseudo-ranks instead of ranks, is 
recommended.

{\color{black} 
A completely different approach to circumventing the above mentioned problem 
could be given by so-called aligned rank procedures (see Hodges and Lehmann, 
1962 or Puri and Sen, 1985). However, these procedures are restricted to 
semi-parametric models and are not applicable to ordinal data since the 
hypotheses are formulated in terms of artificially introduced parameters of 
metric data. In addition, the typical invariance under strictly monotone 
transformations of the data is lost by the two different transformations.

On the contrary, the proposed pseudo-rank procedures share the typical invariance of 
rank-based procedures and do not lead to the paradoxical results described in 
Sections~\ref{owl} and \ref{twl}. They can additionally be used to construct 
confidence intervals for meaningful nonparametric effect sizes as described in 
Section~\ref{expl.confi}. In case of factorial designs with 
independent observations pseudo-rank procedures are already implemented in the {\it R}-package {\it rankFD} by 
choosing the option for {\it unweighted} effects.}

\newpage

\section{References}
\bdes

\item {\sc Akritas, M. G., Arnold, S. F. and Brunner, E. (1997)}. Nonparametric 
hypotheses and rank statistics for unbalanced factorial designs. {\em Journal 
of the American Statistical Association} {\bf 92}, 258-265.  

\item
{\sc Akritas, M.~G. and  Brunner, E. (1997)}. A unified approach to ranks tests 
in mixed models. {\em Journal of Statistical Planning and Inference} {\bf 61},
249--277.

\item {\sc Bewick, V., Cheek, L.}, and {\sc Ball, J.} (2004). Statistics Review 
10: Further Nonparametric Methods. {\em Critical Care} {\bf 8}, 196-199.

\item {\sc Birnbaum, Z.~W. and Klose, O.~M. (1957)}. Bounds for the Variance of 
the Mann-Whitney Statistic. {\em Annals of Mathematical Statistics} {\bf 28}, 
933--945.

\item {\sc Brunner, E., Bathke, A.~C., and Konietschke, F. (2018)}. {\em Rank- 
and Pseudo-Rank Procedures in Factorial Designs - Using R and SAS - Independent 
Observations.} Lecture Notes in Statistics, Springer, Heidelberg.

\item {\sc Brunner, E., Konietschke, F., Pauly, M., and Puri, M.L. (2017)}. 
Rank-Based Procedures in Factorial Designs: Hypotheses about Nonparametric 
Treatment Effects. {\em Journal of the Royal Statistical Society Series B}, 
{\bf 79}, 1463--1485.

\item {\sc Brunner, E. and Munzel, U. (2000).} The Nonparametric Behrens-Fisher 
Problem: Asymptotic Theory and a Small-Sample Approximation. {\em Biometrical 
Journal} {\bf 42}, 17--25.

\item {\sc Brunner, E. and Munzel, U. (2013).} {\em Nichtparametrische 
Datenanalyse.} $2^{nd}$ edition, Springer, Heidelberg.

\item{\sc Brunner, E. and Puri, M.L. (2001)}. Nonparametric Methods in 
Factorial Designs. {\em Statistical Papers} {\bf 42}, 1--52.

\item{\sc Brunner, E. and Puri, M.~L. (2002)}. A class of rank-score tests in 
factorial designs. {\em Journal of Statistical Planning and Inference} 
{\bf 103}, 331--360.

\item {\sc Gao, X. and Alvo, M. (2005$a$)}. A Nonparametric Test for 
Interaction in Two-Way Layouts. {\em The Canadian Journal of Statistics} 
{\bf 33}, 529--543.

\item {\sc Gao, X. and Alvo, M. (2005$b$)}. A Unified Nonparametric Approach 
for Unbalanced Factorial Designs. {\em Journal of the American Statistical 
Association}, {\bf 100}, 926--941.


\item {\sc Govindarajulu, Z. (1968)}. Distribution-free confidence bounds for 
$Pr\{ X < Y \}$. {\em Annals of the Institute of Statistical Mathematics} {\bf 
20}, 229--238.

\item {\sc Hettmansperger, T.~P. and  Norton, R.~M. (1987)}. Tests for 
patterned alternatives in k-sample problems. {\em Journal of the American 
Statistical Association} {\bf 82}, 292--299.

\item {\sc Jonckheere, A.~R. (1954)}. A Distribution-free $k$-sample Test 
Against Ordered Alternatives. {\em Biometrika} {\bf 41}, 133--145.

\item {\sc Konietschke, F. (2009)}. Simultane Konfidenzintervalle f\"ur 
nichtparametrische relative Kontrasteffekte. PhD thesis, Inst. of Math. 
Stochastics, University of G\"ottingen.

\item {\sc Konietschke, F. Hothorn, L.~A., and Brunner, E.} (2012). Rank-based 
multiple test procedures and simultaneous confidence intervals. {\em Electronic 
Journal of Statistics} {\bf 6}, 737--758. DOI: 10.1214/12-EJS691

\item {\sc Kruskal, W.~H. (1952)}. A Nonparametric Test for the Several Sample 
Problem. {\em Annals of Mathematical Statistics} {\bf 23}, 525--540.

\item {\sc Kruskal, W.~H. and  Wallis, W.~A. (1952)}. The use of ranks in 
one-criterion variance analysis. {\em Journal of the American Statistical 
Association} {\bf 47}, 583--621.

\item {\sc Kulle, B. (1999)}. Nichtparametrisches Behrens-Fisher-Problem im 
Mehrstichprobenfall. Diploma Thesis, Inst. of Math. Stochastics, University of 
G\"ottingen.

\item {\sc Peterson, I. (2002)}. Tricky Dice Revisited. {\em Science News} 
{\bf 161}, \\ https://www.sciencenews.org/article/tricky-dice-revisited

\item {\sc Puri, M.~L. (1964)}. Asymptotic efficiency of a class of $c$-sample 
tests. {\em Annals of Mathematical Statistics} {\bf 35}, 102--121.

\item {\sc Puri, M.~L. and Sen, P.~K. (1985)}. {\em Nonparametric methods in 
general linear models}. Wiley, New York.

\item {\sc Ruymgaart, F.~H. (1980)}. A unified approach to the asymptotic 
distribution theory of certain midrank statistics. In: {\em Statistique non 
Parametrique Asymptotique}, 1--18, J.P. Raoult (Ed.), Lecture Notes on 
Mathematics, No. 821, Springer, Berlin.

\item {\sc Sen, P.~K. (1967)}. A note on asymptotically distribution-free 
confidence intervals for $Pr(X<Y )$ based on two independent samples. 
{\em Sankhya, Series A} {\bf 29}, 95--102.

\item {\sc Thangavelu, K. and Brunner, E. (2007).} Wilcoxon-Mann-Whitney test 
for stratified samples and Efron's paradox dice. {\em Journal of Statistical 
Planning and Inference} {\bf 137}, 720--737.

\item {\sc Terpstra, T.~J. (1952)}. The asymptotic normality and consistency of 
Kendall's test against trend, when ties are present in one ranking. 
{\em Indagationes Mathematicae} {\bf 14}, 327--333.

\item {\sc Whitley, E., and Ball, J. (2002).} Statistics Review 6: 
Nonparametric Methods. {\em Critical Care}, {\bf 6(6)}, 509-513.

\item {\sc Zaremba, S.~K. (1962)}. A generalization of Wilcoxon's tests. 
{\em Monatshefte f\"ur Mathematik} {\bf 66}, 359--70.

\edes

\end{document}